\documentstyle[12pt]{amsart}
\makeatletter
\def\cases#1{\left\{\,\vcenter{\normalbaselines\m@th
    \ialign{$##\hfil$&\quad{##}\hfil\crcr#1\crcr}}\right.}
\makeatother

\newtheorem{theorem}{Theorem}\newtheorem{lemma}[theorem]{Lemma}

\let\pont\relax

\def\bs{\bigskip}
\renewcommand{\text}[1]{\quad\mbox{#1}\quad}
\def\beqn{\begin{eqnarray}}\def\eeqn{\end{eqnarray}}
\def\beqno{\begin{eqnarray*}}\def\eeqno{\end{eqnarray*}}
\def\ds{\displaystyle}\def\eps{\varepsilon}\def\ffi{\varphi}
\def\qed{\ifhmode\unskip\nobreak\fi\quad\ifmmode\Box\else$\Box$\fi}
\def\be{\mbox{}\hskip10pt}\def\({\mbox{$($}}\def\){\mbox{$)$}}

\title{An Algorithmic Version of the Blow-up Lemma}

\author{J\'anos Koml\'os}
\address{\hskip-\parindent
J\'anos Koml\'os\\Mathematics Department\\Rutgers University\\
New Brunswick, NJ 08903}
\email{komlos@@math.rutgers.edu}

\author{Gabor N. Sarkozy}
\thanks{Part  of this paper was written while 
Sarkozy was visiting
MSRI Berkeley, as part of the 
Combinatorics Program. Research at MSRI is supported
in part by NSF grant DMS-9022140.}
\address{\hskip-\parindent
Gabor N. Sarkozy\\
Computer Science Department\\
Worcester Polytechnic Institute\\
Worcester, MA 01609}
\email{gsarkozy@@cs.wpi.edu}

\author{Endre Szemer\'edi}
\address{\hskip-\parindent
Endre Szemer\'edi\\
Computer Science Depratment\\
Rutgers University\\New Brunswick, NJ 08903}
\email{szemered@@cs.rutgers.edu}

\begin{document}

\begin{abstract}
Recently we have developed a new method in graph theory
based on the Regularity Lemma.
The method is applied to find certain
spanning subgraphs in dense graphs.
The other main general tool of the method,
beside the Regularity Lemma,
is the so-called Blow-up Lemma (\cite{KSSz-BL}).
This lemma helps to find bounded degree spanning subgraphs
in $\eps$-regular graphs. 
Our original proof of the lemma is not algorithmic,
it applies probabilistic methods.
In this paper we provide an algorithmic version of the Blow-up Lemma.
The desired subgraph, for an $n$-vertex graph,
can be found in time $O(nM(n))$, where $M(n)=O(n^{2.376})$
is the time needed to multiply two $n$ by $n$ matrices
with 0,1 entries over the integers.
We show that the algorithm can be parallelized
and implemented in $NC^5$.
\end{abstract}

\maketitle

\section{Introduction}
\subsection{Notations and definitions}

All graphs are simple, that is, they have no loops or multiple edges.
$v(G)$ is the number of vertices in $G$ (order),
$e(G)$ is the number of edges in $G$ (size).
$deg(v)$ (or $deg_G(v)$) is the degree of vertex $v$
(within the graph $G$), and $deg(v,Y)$ (or $deg_G(v,Y)$)
is the number of neighbours of $v$ in $Y$.
$\delta(G)$ and $\Delta(G)$ are the minimum degree
and the maximum degree of $G$.
$N(x)$ (or $N_G(x)$) is the set of neighbours of the vertex $x$,
and $e(X,Y)$ is the number of edges between $X$ and $Y$.
A bipartite graph $G$ with color-classes $A$ and $B$
and edges $E$ will sometimes be written as $G=(A,B,E)$.
For disjoint $X,Y$, we define the {\bf density}
$$d(X,Y)={e(X,Y)\over|X|\cdot|Y|}\,.$$
The density of a bipartite graph $G=(A,B,E)$ is the number
$$d(G)=d(A,B)={|E|\over|A|\cdot|B|}\,.$$
For two disjoint subsets $A,B$ of $V(G)$,
the bipartite graph with vertex set $A\cup B$ which has
all the edges of $G$ with one endpoint in $A$
and the other in $B$ is called the pair $(A,B)$.

A pair $(A,B)$ is $\eps$-{\bf regular}
if for every $X\subset A$ and $Y\subset B$ satisfying
$$|X|>\eps|A|\text{and}|Y|>\eps|B|$$
we have $$|d(X,Y)-d(A,B)|<\eps.$$
A pair $(A,B)$ is $(\eps,\delta)${\bf-super-regular}
if it is $\eps$-regular
and furthermore,$$deg(a)\geq\delta|B|\text{for all}a\in A,$$
$$\text{and}deg(b)\geq\delta|A|\text{for all}b\in B.$$

$H$ is {\bf embeddable} into $G$ if $G$ has a subgraph isomorphic to $H$,
that is, if there is a one-to-one map (injection)
$\ffi:V(H)\to V(G)$ such that $\{x,y\}\in E(H)$ implies
$\{\ffi(x),\ffi(y)\}\in E(G)$.

As the model of computation we choose the weakest possible version of
a {\em PRAM}, in which concurrent reads or writes of the same location are
not allowed ({\em EREW}, see \cite{KR} for a discussion of the various
{\em PRAM}
models.) When researchers investigate the parallel complexity of a
problem, the main question is whether a polylogarithmic running time
is achievable on a {\em PRAM} containing a polynomial number of processors.
If the answer is positive then the problem and the corresponding
algorithm are said to belong to class $NC$ introduced in \cite{PI} (see also
\cite{CO},\cite{VA}). When the running time is $O((\log n)^i)$, the
algorithm is in $NC^i$.

\subsection{An Algorithmic Version of the Blow-up Lemma}

In recent years the interaction between combinatorics and the theory
of algorithms is getting stronger and stronger. It is therefore not
surprising that there has been a significant interest in
converting existence proofs into efficient algorithms.
Many examples of this type can be found in \cite{A1, PA}.
Some of these are general methods, so algorithmic versions of these methods
immediately imply efficient algorithms for several problems.
One example is the Lov\'asz Local Lemma whose algorithmic aspects have been
studied recently in \cite{A2, B}.
Another example is the Regularity Lemma \cite{Sz1}.
The basic content of this lemma could be described
by saying that every graph can, in some sense,
be well approximated by random graphs.
Since random graphs of a given edge density
are much easier to treat than all graphs
of the same edge-density, the Regularity Lemma helps us
to carry over results that are trivial for random graphs
to the class of all graphs with a given number of edges.
The lemma has numerous applications in various areas
including combinatorial number theory \cite{FGR,SZ2},
computational complexity \cite{HMT} and extremal graph theory
\cite{AY1,BESS,CRST,EFR,FR,FF,F,FGK,R,SiSo}.
Recently an $NC^1$-algorithmic version was given in \cite{ADLRY}.

During the past couple of years we have developed
a new method in graph theory based on the Regularity Lemma.
We usually apply this method to find
certain spanning subgraphs in dense graphs.
Typical examples are spanning trees 
(Bollob\'as-conjecture, see \cite{KSSz-B}),
Hamiltonian cycles or powers of Hamiltonian cycles
(P\'osa-Seymour conjecture, see \cite{KSSz-PS1,KSSz-PS2})
or $H$-factors for a fixed graph $H$ 
(Alon-Yuster conjecture, see \cite{KSSz-AY}).
The other main general tool in the method,
beside the Regularity Lemma,
is the so-called Blow-up Lemma (\cite{KSSz-BL}).
This lemma helps to find bounded degree spanning subgraphs
in $\eps$-regular graphs. 
The rough idea of the original proof of this lemma
was the following: we use a randomized greedy algorithm
to embed most of the vertices of the bounded degree graph,
and then finish the embedding by a K\"onig-Hall argument.
Given the recent algorithmic version of the Regularity Lemma,
the obvious question arises whether there is an algorithmic version
also for the Blow-up Lemma.
In this paper we give an affirmative answer to this question.
\begin{theorem}{\bf (An Algorithmic Version of the Blow-up Lemma)}
\label{blowup}\pont
Given a graph $R$ of order $r$ and positive parameters $\delta,\Delta$,
there exists an $\eps>0$ such that the following holds.
Let $N$ be an arbitrary positive integer,
and let us replace the vertices of $R$
with pairwise disjoint $N$-sets
$V_1,V_2,\ldots,V_r$ \(blowing up\).
We construct two graphs on the same vertex-set $V=\cup V_i$.
The graph $R(N)$ is obtained by replacing all edges of $R$
with copies of the complete bipartite graph $K_{N,N}$,
and a sparser graph $G$ is constructed by replacing
the edges of $R$ with some $(\eps,\delta)$-super-regular pairs.
If a graph $H$ with $\Delta(H)\leq\Delta$ is embeddable into $R(N)$
then it is already embeddable into $G$.
We can construct a copy of $H$ in $G$ in $O(nM(n))$ sequential time,
where $M(n)=O(n^{2.376})$ is the time needed
to multiply two $n$ by $n$ matrices
with 0,1 entries over the integers.
Furthermore, the algorithm can be parallelized
and implemented in $NC^5$.
\end{theorem}
{\bf Remark.}
% The theorem is also true for unequal vertex-clusters $V_i$.
For some very special cases of this theorem
(e.g. a Hamiltonian path in a super-regular pair)
$NC$ algorithms can be found in \cite{S1}.

When using the Blow-up Lemma, we typically also need
the following strengthened version:
Given $c>0$, there are positive functions
$\eps=\eps(\delta,\Delta,r,c)$ and
$\alpha=\alpha(\delta,\Delta,r,c)$
such that the Blow-up Lemma remains true
if for every $i$ there are certain vertices $x$
to be embedded into $V_i$ whose images are {\it a priori}
restricted to certain sets $C_x\subset V_i$ provided that\newline
\be (i) each $C_x$ within a $V_i$ is of size at least $c|V_i|$,\newline
\be (ii) the number of such restrictions within a $V_i$
is not more than $\alpha|V_i|$.

Our proof is going to be similar to our original probabilistic proof,
but we have to replace the probabilistic arguments with deterministic ones.
In Section 2 we give a deterministic sequential algorithm
for the embedding problem without considering implementation
and time complexity issues.
In Section 3 we show that the algorithm is correct.
Implementation is discussed in Section 4.
Finally, Section 5 contains various algorithmic applications.

\section{The algorithm}
The main idea of the algorithm is the following.
We embed the vertices of $H$ one-by-one by following
a greedy algorithm, which works smoothly
until there is only a small proportion of $H$ left,
and then it may get stuck hopelessly.
To avoid that, we will set aside a positive proportion
of the vertices of $H$ as buffer vertices.
Most of these buffer vertices will be embedded
only at the very end by using a K\"onig-Hall argument.

\subsection{Preprocessing}
We will assume that $|V(H)|=|V(G)|=|\cup_iV_i|=n=rN$.
We can also assume, without restricting generality,
that $N\geq N_0(\delta,\Delta,r)$,
for Theorem \ref{blowup} is trivial for small $N$
since $\eps$-regularity with a small enough $\eps$ implies $G=R(N)$.
Finally, we will assume for simplicity, that the density
of every super-regular pair in $G$ is exactly $\delta$.
This is not a significant restriction,
otherwise we just have to put everywhere the actual density
instead of $\delta$.

We will use the following parameters:
$$\eps\ll\eps'\ll\eps''\ll\delta'''\ll\delta''\ll\delta'\ll\delta,$$
where $a\ll b$ means that $a$ is small enough compared to $b$.

For easier reading, we will mostly use the letter $x$
for vertices of $H$,
and the letter $v$ for vertices of the host graph $G$.

Given an embedding of $H$ into $R(N)$, it defines an {\it assignment}
$$\psi:V(H)\to\{V_1,V_2,\ldots,V_r\},$$
and we want to find an {\it embedding}
$$\ffi:V(H)\to V(G),\quad\ffi\ \ \mbox{is one-to-one}$$
such that $\ffi(x)\in\psi(x)$ for all $x\in V(H)$.
We will write $X_i=\psi^{-1}(V_i)$ for $i=1,2,\dots,r$.

Before we start the algorithm, we order the vertices
of $H$ into a sequence $S=(x_1,x_2,\ldots,x_n)$
which is more or less, but not exactly,
the order in which the vertices will be embedded.
Let $m=r\lceil\delta'N\rceil$.
For each $i$, choose a set $B_i$
of $\lceil\delta'N\rceil$ vertices in $X_i$ such that
any two of these vertices
are at a distance at least four in $H$.
(This is possible, for $H$ is a bounded degree graph.)
These vertices $b_1,\dots,b_m$ will be called
the {\bf buffer vertices}
and they will be the last vertices in $S$.

The order $S$ starts with the neighbourhoods
$N_H(b_1),N_H(b_2),\dots,N_H(b_m)$.
The length of this initial segment of $S$
will be denoted by $T_0$.
Thus $T_0=\ds\sum_{i=1}^m|N_H(b_i)|\leq\Delta m$.

The rest of $S$ is an arbitrary ordering
of the leftover vertices of $H$.

(When certain images are {\it a priori} restricted,
as in the remark after the theorem,
we also list the restricted vertices at the beginning of $S$
right after the neighbours of the buffer vertices.)

\subsection{Sketch of the algorithm}
In {\em Phase 1} of the algorithm we will embed
the vertices in $S$ one-by-one into $G$
until all non-buffer vertices are embedded.
For each $x_j$ not embedded yet
(including the buffer vertices)
we keep track of an ever shrinking host set $H_{t,x_j}$
that $x_j$ is confined to at time $t$,
and we only make a final choice for the location of $x_j$
from $H_{t,x_j}$ at time $j$.
At time 0, $H_{0,x_j}$ is the cluster
that $x_j$ is assigned to.
For technical reasons we will also maintain another similar set,
$C_{t,x_j}$, where we will ignore the
possibility that some vertices are occupied already.

In {\em Phase 2}, we embed the leftover vertices
by using a K\"onig-Hall type argument.

\subsection{Embedding Algorithm}

At time 0, set $C_{0,x}=H_{0,x}=\psi(x)$ for all $x\in V(H)$.
Put $T_1=\lfloor\delta''n\rfloor$.

\bs{\bf Phase 1}\\
For $t\geq 1$, repeat the following steps.

{\em Step 1} - {\em Extending the embedding.}\quad
We embed $x_t$. Consider the vertices in $H_{t-1,x_t}$.
We will pick one of these vertices as the image $\ffi(x_t)$
by using the Selection Algorithm
(described below in Section \ref{selection}).

{\it Step 2} - {\em Updating.}\quad
For each unembedded vertex $y$
(i.e. the set of vertices $x_j,t<j\leq n$), set
$$C_{t,y}=\cases{C_{t-1,y}\cap N_G(\ffi(x_t))
&if $\{x_t,y\}\in E(H)$\cr
C_{t-1,y}&otherwise,}$$and
$$H_{t,y}=\cases{H_{t-1,y}
\cap N_G(\ffi(x_t))&if $\{x_t,y\}\in E(H)$\cr
H_{t-1,y}\setminus\{\ffi(x)\}&otherwise.}$$

{\em Step 3} - {\em Exceptional vertices in $H$.}\quad

1. If $T_1$ does not divide $t$, then go to Step 4.

2. If $T_1$ divides $t$, then we do the following.
We find all exceptional unembedded vertices $y\in H$ such that
$|H_{t,y}|\leq(\delta')^2n$.
At this point we slightly change the order
of the remaining unembedded vertices in $S$.
We bring these exceptional vertices forward
(even if they are buffer vertices), followed by
the non-exceptional vertices in the same relative order as before.
For simplicity we still use the notation $(x_1,x_2,\ldots,x_n)$
for the new order.

{\em Step 4} - {\em Exceptional vertices in $G$.}\quad

1. If $t\not=T_0$, then go to Step 5.

2. If $t=T_0$, then we do the following.
We find the set (denoted by $E_i$)
of those exceptional vertices $v\in V_i$, $1\leq i\leq r$
for which $v$ is not covered yet in the embedding and
$$\left|\left\{b:b\in B_i,v\in C_{t,b}\right\}\right|
<\delta''|B_i|.$$
Once again we are going to change slightly the order
of the remaining unembedded vertices in $S$.
We choose a set $E$ of vertices
$x\in H$ of size $\sum_{i=1}^{r}|E_i|$
(more precisely $|E_i|$ vertices from $X_i$
for all $1\leq i\leq r$) with
$$H_{t,x}=H_{0,x}\setminus\{\ffi(x_j):j\leq t\}
=\psi(x)\setminus\{\ffi(x_j):j\leq t\}.$$
Thus in particular, if $x\in X_i$, then $E_i\subset H_{t,x}$.
For example, we may choose the vertices in $E$
as vertices in $H$ that are at a distance at least four from each other
and any of the vertices embedded so far.
We are going to show later in the proof of correctness
that this is possible since $H$ is a bounded degree graph
and $\sum_{i=1}^r|E_i|$ is very small.
We bring the vertices in $E$ forward,
followed by the remaining unembedded vertices
in the same relative order as before.
Again, for simplicity we keep the notation
$(x_1,x_2,\ldots,x_n)$ for the resulting order.
After this we set $t\leftarrow t+1$ and go back to Step 1.
However, when we perform Step 1 for a vertex $x\in E$
we pick $\ffi(x)$ from the appropriate $E_i$.

{\em Step 5} -
If there are no more unembedded non-buffer vertices left,
then set $T=t$ and go to Phase 2,
otherwise set $t\leftarrow t+1$ and go back to Step 1.

\bs{\bf Phase 2}\\
Find a system  of distinct representatives of the sets
$H_{T,y}$ for all unembedded $y$
(i.e. the set of vertices $x_j$, $T<j\leq n$).

\subsection{Selection Algorithm}\label{selection}

We choose a vertex $v\in H_{t-1,x_t}$ as the image $\ffi(x_t)$
for which the following hold for all unembedded $y$
with $\{x_t,y\}\in E(H)$
$$(\delta-\eps)|H_{t-1,y}|\leq deg_G(v,H_{t-1,y})
\leq(\delta+\eps)|H_{t-1,y}|,$$
$$(\delta-\eps)|C_{t-1,y}|\leq deg_G(v,C_{t-1,y})
\leq(\delta+\eps)|C_{t-1,y}|$$
and
$$(\delta-\eps)|C_{t-1,y}\cap C_{t-1,y'}|
\leq deg_G(v,C_{t-1,y}\cap C_{t-1,y'})
\leq(\delta+\eps)|C_{t-1,y}\cap C_{t-1,y'}|,$$
for at least $(1-\eps')$ proportion
of the unembedded vertices $y'$ with $\psi(y')=\psi(y)$.

\section{Proof of correctness}

The following claims state that our algorithm finds
a good embedding of $P$ into $G$.

\bigskip\noindent
{\bf Claim 1.}\ {\em Phase 1 always succeeds.}

\bigskip\noindent
{\bf Claim 2.}\ {\em Phase 2 always succeeds.}

\bigskip\noindent
If at time $t$,
$S$ is a set of unembedded vertices $x\in H$ with $\psi(x)=V_i$,
then we define the bipartite graph $U_t$ as follows.
One color class is $S$, the other is $V_i$, and we have an edge
between a $x\in S$ and a $v\in V_i$ whenever $v\in C_{t,x}$.

In the proofs of the above claims the following lemma
will play a major role.

\begin{lemma}\label{egy}\pont
We are given integers $1\leq i\leq r$, $1\leq t\leq T$
and a set $S\subset X_i$ of unembedded vertices at time $t$ with
$|S|\geq(\delta''')^2|X_i|=(\delta''')^2N$.
Then apart from an exceptional set $F$ of size at most $\eps''N$,
for every vertex $v\in V_i$ we have the following
% $$(1+\eps'')d(U_t)|S|\geq 
$$deg_{U_t}(v)
=\left|\{x:x\in S,v\in C_{t,x}\}\right|\geq(1-\eps'')d(U_t)|S|
\quad\left(\geq\frac{\delta^{\Delta}}2|S|\right).$$
\end{lemma}

{\bf Proof.} In the proof of this lemma we will use the
``defect form'' of the Cauchy-Schwarz inequality
(just as in the original proof of the Regularity Lemma): if
$$\sum_{k=1}^m X_k=\frac{m}{n}\sum_{k=1}^n X_k
+\Delta\qquad(m\leq n)$$
then
$$\sum_{k=1}^n X_k^2\geq\frac1n
\left(\sum_{k=1}^nX_k\right)^2+\frac{\Delta^2n}{m(n-m)}.$$

Assume indirectly that the statement in Lemma \ref{egy} is not true,
that is, $|F|>\eps''N$.
Let us write $\nu(t,x)$ for the number of neighbors (in $H$)
of $x$ embedded by time $t$.
Then in $U_t$ we have the following.
\begin{equation}\label{nagy}e(U_t)=d(U_t)|S||V_i|
=\sum_{v\in V_i}deg_{U_t}(v)=\sum_{x\in S}deg_{U_t}(x)
=\sum_{x\in S}|C_{t,x}|\geq\sum_{x\in S}(\delta-\eps)^{\nu(t,x)}N.\end{equation}
We also have
\begin{equation}\label{kicsi}
\sum_{x_1\in S}\sum_{x_2\in S}|N_{U_t}(x_1)\cap N_{U_t}(x_2)|
\leq\sum_{x_1\in S}\sum_{x_2\in S}
(\delta+\eps)^{\nu(t,x_1)+\nu(t,x_2)}N+4\Delta\eps'N|S|^2.
\end{equation}

On the other hand using (\ref{nagy}) and the Cauchy-Schwarz
inequality with $m=\lceil\eps''N\rceil$, we get
$$\sum_{x_1\in S}\sum_{x_2\in S}|N_{U_t}(x_1)\cap N_{U_t}(x_2)|
=\sum_{v\in V_i}(deg_{U_t}(v))^2\geq$$
$$\geq\frac1N\left(\sum_{v\in V_i}deg_{U_t}(v)\right)^2
+(\eps'')^3d(U_t)^2N|S|^2\geq$$
$$\geq\frac1N\left(\sum_{x\in S}
(\delta-\eps)^{\nu(t,x)}N\right)^2
+(\eps'')^3\delta^{2\Delta}N|S|^2\geq$$
$$\geq\sum_{x_1\in S}\sum_{x_2\in S}
(\delta-\eps)^{\nu(t,x_1)+\nu(t,x_2)}N
+(\eps'')^3\delta^{2\Delta}N|S|^2,$$
which is a contradiction with (\ref{kicsi}), if $\eps''$
is sufficiently large compared to $\eps'$ and $\eps$.

An easy consequence of Lemma \ref{egy} is the following lemma.

\begin{lemma}\label{ketto}\pont
We are given integers $1\leq i\leq r$, $1\leq t\leq T$,
a set $S\subset X_i$ of unembedded vertices at time $t$ with
$|S|\geq\delta'''|X_i|=\delta'''N$ and a set $A\subset V_i$ with
$|A|\geq \delta'''|V_i|=\delta'''N$. Then
apart from an exceptional set $F$ of size at most $(\delta''')^2N$,
for every vertex $x\in S$ we have the following
\begin{equation}\label{jolmetsz}
|A\cap C_{t,x}|\geq\frac{|A|}{2N}|C_{t,x}|.
\end{equation}
\end{lemma}

{\bf Proof.} Assume indirectly that the statement is not true, i.e.
there exists a set $S'\subset S$ with $|S'|>(\delta''')^2N$
such that for every $x\in S'$ (\ref{jolmetsz}) does not hold.
Once again we consider the bipartite graph $U_t=U_t(S',V_i)$.
We have
$$\sum_{v\in A}deg_{U_t}(v)=\sum_{x\in S'}
|A\cap C_{t,x}|<\frac{|A|}{2N}\,\sum_{x\in S'}|C_{t,x}|
=\frac{|A|}{2N}\ d(U_t)|S'|N.$$
On the other hand, applying Lemma \ref{egy} for $S'$ we get
$$\sum_{v\in A}deg_{U_t}(v)\geq(1-\eps'')d(U_t)|S'|(|A|-\eps''N)$$
contradicting the previous inequality.

Finally we have
\begin{lemma}\label{harom}\pont
For every $1\leq t\leq T$
and for every vertex $y$ that is unembedded at time $t$,
we have$$|H_{t,y}|>\delta''N.$$
\end{lemma}

{\bf Proof.} Indeed, Lemma \ref{ketto}
implies that every time option 2 is executed
in Step 3 (i.e. $T_1|t$), we find at most
$(\delta''')^2N\ (\ll\delta''N)$ exceptional vertices.
This fact clearly implies Lemma \ref{harom}.

\bigskip
Lemma \ref{egy} implies that for $E_i$ in Step 3 of Phase 1
of the Embedding Algorithm we have $|E_i|<\eps''N$, thus
$$\sum_{i=1}^{r}|E_i|\leq r\eps''N.$$

Lemma \ref{harom} implies that the Selection Algorithm
always selects from a set
of size at least $\delta''N$. Furthermore, since
$|H_{t,y}|>\eps N$ for every $1\leq t\leq T$ and 
for every vertex $y$ unembedded at time $t$,
the $\eps$-regularity and some simple computation imply
that at most $\eps'N$
($\ll\delta''N$ if $\eps'$ is small enough)
vertices violate the degree requirements
in the Selection Algorithm.
Therefore the Selection Algorithm always succeeds
in selecting an image $\ffi(x_t)$, proving Claim 1.

\bigskip\noindent
{\bf Proof of Claim 2.}
We want to show that we can find a system of distinct
representatives of the sets $H_{T,x_j},T<j\leq n$,
where the $H_{T,x_j}$-s belong to a given cluster $V_i$.

To simplify notation, let us denote by $Y$
the set of remaining vertices in $V_i$,
and by $X$ the set of remaining
unembedded (buffer) vertices assigned to $V_i$.
If $x=x_j\in X$ then write $H_x$ for its possible location
$H_{T,x_j}$ at time $T$. Also write $M=|X|=|Y|$.

The K\"onig-Hall condition for the existence of a system of
distinct representatives obviously follows from the following
three conditions:
\begin{equation}\label{konig1}
|H_x|>\delta'''M\quad\mbox{for all}\quad x\in X,
\end{equation}
\begin{equation}\label{konig2}
|\displaystyle\mathop\cup_{x\in S}H_x|\geq(1-\delta''')M
\text{for all subsets}S\subset X,|S|\geq\delta'''M,
\end{equation}
\begin{equation}\label{konig3}
|\displaystyle\mathop\cup_{x\in S}H_x|=M
\text{for all subsets}S\subset X,|S|\geq(1-\delta''')M.
\end{equation}
(\ref{konig1}) is an immediate consequence of Lemma \ref{harom},
(\ref{konig2}) is a consequence of Lemma \ref{egy}.

Finally to prove (\ref{konig3}), we have to show that every vertex
in $Y\subset V_i$ belongs to at least $\delta''' |X|$ sets $H_x$.
But this is trivial from the construction of the embedding algorithm,
in Step 3 of Phase 1 we took care
of the small number of exceptional vertices
for which this is not true.
This finishes the proof of Claim 2 and the proof of correctness.

\section{Implementation}

The sequential implementation is immediate.
In Phase 1 we have $\leq n$ iterations,
and it is not hard to see that one iteration
can be implemented in $O(M(n))$ time.
Phase 2 can be implemented in
$O(rN^{5/2})=O(M(n))$ time by applying an algorithm
for finding a maximum matching in a bipartite graph
(see e.g. \cite{HK,LP}).

For the parallel implementation, our main tool is the $NC^4$ algorithm
for the maximal independent set problem.
A subset $I$ of the vertices of a graph $G$ is independent
if there are no edges between any two vertices in $I$.
An independent set $I$ is {\em maximal}
if it is not a proper subset of any other independent set.
Karp and Wigderson (\cite{KW}) were the first to give an
$NC^4$-algorithm for this problem.
Better algorithms were later described in \cite{ABI,GS1,GS2,L}.
We call this the MIS algorithm.

For the parallelization of the embedding algorithm, we show that if
$\alpha$ is a small enough constant and $n'$ is the number of
remaining unembedded vertices,
then we can embed $\lfloor\alpha n'\rfloor$ vertices
in parallel.
First we pick these vertices by running the MIS algorithm
on the following auxiliary graph.
The vertices are the vertices of $H$,
and we put an edge between two vertices,
if either they are at a distance less than 4,
or both vertices are embedded already.
If in the maximal independent set that we find,
we have a vertex that is embedded already
(we can have only one such vertex),
then we remove this vertex from the independent set.
We keep $\lfloor\alpha n'\rfloor$ vertices from the remaining vertices.
These vertices are brought forward in the order in the embedding algorithm
and we embed these vertices in parallel. For each such vertex we determine
the set of vertices where it could be embedded by the embedding algorithm.
Once again running the MIS algorithm on the appropriately defined
auxiliary graph, we can choose a distinct representative from these sets.
Finally we embed each vertex to its representative.
We iterate this procedure until the number
of remaining unembedded vertices is $\leq(\log n)^5$,
and then we embed these vertices sequentially.
Thus Phase 1 can be implemented in
$O(\log n(\log n)^4)=O((\log n)^5)$ parallel time.

For Phase 2 it remains to show, that if the bipartite graph $U_t$
is defined as above between $X$ and $Y$
(i.e. there is an edge between $x\in X$ and $v\in Y$
if and only if $v\in H_x$), then we can construct
a perfect matching in $U_t$ in $NC^4$. For this purpose we obtain
a maximal matching by running MIS on the linegraph of $U_t$.
Then obviously for the remaining unmatched vertices,
say $Z(X)$ and $Z(Y)$, $|Z(X)|=|Z(Y)|$, and $Z(X)\cup Z(Y)$
is an independent set. From (\ref{konig2})
$|Z(X)|\leq\delta'''N$ follows.
Furthermore, Lemma \ref{ketto} and Lemma \ref{harom} imply
that, if we take $x\in Z(X), v\in Z(Y)$,
then there are many ($\gg\delta'''N$) (internally) vertex-disjoint
alternating paths of length 5 between $x$ and $v$.
Running again MIS on an appropriately defined auxiliary graph,
we can find vertex-disjoint alternating paths of length 5
between the pairs in $Z(X)\cup Z(Y)$.
Changing the matching edges to non-matching edges on these paths
we get a perfect matching.

\section{Applications}

In most applications of our method, the only non-constructive parts
are the Regularity Lemma and the Blow-up Lemma.
Therefore, the existence proofs together with the $NC^1$ version
of the Regularity Lemma and Theorem \ref{blowup},
provide several immediate algorithmic applications.
Let us mention here three applications.
Additional applications and the details of the proofs
will appear in the full version of the paper.

\begin{theorem}\label{bollobas}\pont
(Existential version in \cite{KSSz-B},
NC-version in \cite{S1})
Let $\Delta$ and $0<\delta<1/2$ be given.
Then there exists an $n_0$ with the following properties.
If $n\geq n_0$, $T$ is a tree of order $n$ and maximum degree
$\Delta$, and $G$ is a graph of order $n$ and minimum degree
at least $((1/2)+\delta)n$, then $T$ is a subgraph of $G$.
Furthermore, a copy of $T$ in $G$ can be found
in $O(nM(n))$ sequential time as well as in $NC^5$.
\end{theorem}

\begin{theorem}\label{posa}\pont
(Existential version in \cite{KSSz-PS2})
There exists an $n_0$ such that
if $G$ has order $n$ with $n\geq n_0$ and
$\delta(G)\geq2n/3,$
then $G$ contains the square of a Hamiltonian cycle.
Furthermore, a copy of a square of a Hamiltonian cycle
can be found in $O(nM(n))$ sequential time as well as in $NC^5$.
\end{theorem}

\begin{theorem}\label{seymour}\pont
(Existential version in \cite{KSSz-PS1})
For any $p>0$ and positive integer $k$
there exists an $n_0=n_0(p,k)$
such that if $G$ has order $n\geq n_0$ and minimal degree
$$\delta(G)\geq\left(\frac{k}{k+1}+p\right)n,$$
then $G$ contains the $k^{th}$ power of a Hamiltonian cycle.
Furthermore, a copy of the $k^{th}$ power of a Hamiltonian cycle
can be found in $O(nM(n))$ sequential time as well as in $NC^5$.
\end{theorem}


\begin{thebibliography}{99}
\bibitem{A1}N. Alon, {\em Non-constructive proofs in combinatorics}, 
in: Proceedings of the International Congress of Mathematicians, Kyoto 1990,
Japan, 1421-1429, Springer Verlag, Tokyo, 1991.

\bibitem{A2}N. Alon, {\em A parallel algorithmic version of the local lemma},
Random Struct. Algorithms 2 (1991), 367-378,
in: Proc. 32nd IEEE FOCS (1991), 586-593.

\bibitem{ADLRY}N. Alon, R. Duke, H. Leffman, V. R\"odl, R. Yuster,
{\em The algorithmic aspects of the regularity lemma},
Proc. 33rd IEEE FOCS (1993), 479-481,
Journal of Algorithms 16 (1994), 80-109.

\bibitem{ABI}N. Alon, L. Babai, A. Itai,
{\em A fast randomized parallel algorithm
for the maximal independent set problem},
J. Algorithms 7 (1986), 567-583.

\bibitem{AY1}N. Alon, R. Yuster,
{\em Almost $H$-factors in dense graphs},
Graphs and Combinatorics 8 (1992), 95-102.

\bibitem{AY2}N. Alon, R. Yuster,
{\em $H$-factors in dense graphs},
J. Combinatorial Theory Ser. B, to appear.

\bibitem{B}J. Beck,
{\em An algorithmic approach to the Lov\'asz Local Lemma},
Random Struct. Algorithms 2 (1991), 343-365.

\bibitem{BESS}B. Bollob\'as, P. Erd\H{o}s,
M. Simonovits, E. Szemer\'edi,
{\em Extremal graphs without large forbidden subgraphs},
in: Advances in Graph Theory, Cambridge Combinatorial Conference,
Trinity College, 1977, Ann. Discrete Math. 3 (1978), 29-41.

\bibitem{CRST}V. Chv\'atal, V. R\"odl, E. Szemer\'edi, W. T. Trotter Jr.,
{\em The Ramsey number of a graph with bounded maximum degree},
Journal of Combinatorial Theory B34 (1983), 239-243.

\bibitem{CO}S. A. Cook,
{\em Taxonomy of problems with fast parallel algorithms},
Information and Control 64 (1985), 2-22.

\bibitem{EFR}P. Erd\H{o}s, P. Frankl, V. R\"odl, 
{\em The asymptotic number of graphs not containing a fixed subgraph
and a problem for hypergraphs having no exponent}, 
Graphs Combin. 2 (1986), 113-121.

\bibitem{FR}F. Franek, V. R\"odl,
{\em On Erd\H{o}s conjecture for multiplicities of complete subgraphs},
Graphs Combin., to appear.

\bibitem{FGR}P. Frankl, R. L. Graham, V. R\"odl,
{\em On subsets of Abelian groups with no 3-term arithmetic progression},
J. Combin. Theory A45 (1987), 157-161.

\bibitem{FF}P. Frankl, Z. F\"uredi,
{\em Exact solution of some Tur\'an-type problems},
J. Combin. Theory A45 (1987), 226-262.

\bibitem{F}Z. F\"uredi,
{\em The maximum number of edges in a minimal graph of diameter 2},
J. Graph Theory 16 (1992), 81-98.

\bibitem{FGK}Z. F\"uredi, M. Goemans, D. J. Kleitman,
{\em On the maximum number of triangles in wheel-free graphs}, to appear.

\bibitem{GS1}M. Goldberg, T. Spencer,
{\em A new parallel algorithm for the maximal independent set problem},
SIAM J. Computing 18 (1989), 419-427.

\bibitem{GS2}M. Goldberg, T. Spencer,
{\em Constructing a maximal independent set in parallel},
SIAM J. Disc. Math. 2 (1989), 322-328.

\bibitem{HMT}A. Hajnal, W. Maass, G. Tur\'an,
{\em On the communication complexity of graph properties},
in: Proc. 20th ACM STOC (1988), 186-191.

\bibitem{HK}J. E. Hopcroft, R. M. Karp,
{\em An $n^{5/2}$ algorithm for maximum matchings
in bipartite graphs}, SIAM J. Comput. 2 (1973), 225-231.

\bibitem{KR}R. Karp, V. Ramachandran,
{\em Parallel algorithms for shared memory machines},
in: Handbook of Theoretical Computer Science,
(J. Van Leeuven ed.) North Holland (1990), 869-941.

\bibitem{KW}R. Karp, A. Wigderson,
{\em A fast parallel algorithm for the
maximal independent set problem}, JACM 32 (1985), 762-773.

\bibitem{KSSz-B}J. Koml\'os, G. N. S\'ark\"ozy, E. Szemer\'edi,
{\em Proof of a packing conjecture of Bollob\'as},
Combinatorics, Probability and Computing 4 (1995), 241-255.

\bibitem{KSSz-BL}J. Koml\'os, G. N. S\'ark\"ozy, E. Szemer\'edi,
{\em The Blow-up Lemma}, to appear in Combinatorica.

\bibitem{KSSz-PS1}J. Koml\'os, G. N. S\'ark\"ozy, E. Szemer\'edi,
{\em On the P\'osa-Seymour conjecture},
to appear in the Journal of Graph Theory.

\bibitem{KSSz-PS2}J. Koml\'os, G. N. S\'ark\"ozy, E. Szemer\'edi,
{\em On the square of a Hamiltonian cycle in dense graphs},
Random Structures and Algorithms 9 (1996), 193-211.

\bibitem{KSSz-AY}J. Koml\'os, G. N. S\'ark\"ozy, E. Szemer\'edi,
{\em Proof of the Alon-Yuster conjecture}, in preparation.

\bibitem{KomSim}J. Koml\'os, M. Simonovits,
{\em Szemer\'edi's Regularity Lemma and its applications in graph theory},
in: Paul Erd\H{o}s is 80, Vol.2.,
Proc. Colloq. Math. Soc. J\'anos Bolyai (1996), 295-352.

\bibitem{LP}L. Lov\'asz, M. D. Plummer, {\em Matching Theory},
(Annals Discr. Math. 29), North Holland, 1986.

\bibitem{L}M. Luby, {\em A simple parallel algorithm
for the maximal independent set problem},
SIAM J. Computing 15 (1986), 1036-1053.

\bibitem{PA}C. H. Papadimitriou,
{\em On graph theoretic lemmata and complexity classes},
Proc. 31st IEEE FOCS (1990), 794-801,
{\em On the complexity of the parity argument
and other inefficient proofs of existence},
J. Comput. System Sci. 47 (1993).

\bibitem{PI}N. J. Pippenger,
{\em On simultaneous resource bounds}, 
Proc. 20th IEEE FOCS (1979), 307-311.

\bibitem{R}V. R\"{o}dl,
{\em On universality of graphs with uniformly distributed edges},
Disc. Math. 59 (1986), 125-134.

\bibitem{S1}G. N. S\'ark\"ozy,
{\em Fast parallel algorithms for finding
Hamiltonian cycles and trees in graphs},
Technical Report 93-81, DIMACS, Rutgers University.

\bibitem{SiSo}M. Simonovits, V. T. S\'os,
{\em Szemer\'edi's partition and quasirandomness},
Random Sruct. Algorithms 2 (1991), 1-10.

\bibitem{Sz1}E. Szemer\'edi, {\em Regular partitions of graphs},
in: Colloques Internationaux C.N.R.S. (J.-C. Bermond, J.-C. Fournier,
M. Las Vergnas, D. Sotteau, Eds.), (1978), 399-401.

\bibitem{SZ2}E. Szemer\'edi,
{\em On a set containing no $k$ elements in arithmetic progression},
Acta Arithmetica XXVII (1975), 199-245. 

\bibitem{VA}L. G. Valiant,
{\em Parallel computation},
Proc. 7th IBM Symp. on Math. Found. of Comp. Science, 1982.
\end{thebibliography}
\end{document}